\documentclass[12pt]{article}
\usepackage{graphicx}
\usepackage{amssymb}
\usepackage{amsmath}
\newcommand{\norm}[1]{\left\Vert#1\right\Vert}

\numberwithin{equation}{section}
\def\be{\begin{equation}\label}
\def\e{\end{equation}}

\author{Lev Sakhnovich}

\title{Infinite Hankel Block Matrices, Extremal  Problems}

\date{}

\begin{document}

\maketitle

Address:
735 Crawford ave., Brooklyn, 11223, New York, USA.\\
 E-mail address: lev.sakhnovich@verizon.net
\begin{center}Abstract \end{center}
In this paper we use the matrix analogue of eigenvalue $\rho_{min}^{2}$ to formulate and to solve
the extremal Nehary problem. 
When $\rho_{min}$ is a scalar, our approach coincides
with  Adamjan-Arov-Krein approach.
 
 \section{Introduction}In the paper we consider a matrix version of the extremal
 Nehary problem [1],[4]. Our approach is based on the notion of a  matrix analogue of the eigenvalue $\rho_{min}^{2}$. The notion of  $\rho_{min}^{2}$  was used in a number
 of the extremal  interpolation problems [2],[3],[7]. We  note that  $\rho_{min}^{2}$ is a solution of a
 non-linear matrix inequality of the Riccati type [2],[6], [7]. Our approach coincides with the 
 Adamjan-Arov-Krein approach [1],
when $\rho_{min}^{2}$ is a scalar matrix.\\ Now we  introduce the main definitions.
 Let $H$ be a fixed separable Hilbert space.
 By ${\ell}_{2}(H)$ we denote the Hilbert space of the sequences
 $\xi=\{\xi_{k}\}_{1}^{\infty},$ where $\xi_{k}{\in}H$ and
 \be{1}
 {\norm{\xi}}^{2}=\sum_{k=1}^{\infty}{\norm{\xi_{k}}}^{2}<\infty.
 \e
 The space of the bounded linear operators acting from  ${\ell}_{2}(H_{1})$ into  ${\ell}_{2}(H_{2})$
 is denoted by  $[{\ell}_{2}(H_{1}), {\ell}_{2}(H_{2})].$ The Hankel operator
 $\Gamma{\in}[{\ell}_{2}(H_{1}), {\ell}_{2}(H_{2})]$ has the form
 \be{2}\Gamma=\{\gamma_{j+k-1}\},\quad 1{\leq}j,k{\leq}\infty,\quad \gamma_{k}{\in}[H_{1},H_{2}].\e
 Let $L_{\infty}[H_{1},H_{2}]$ be the space of the measurable operator -valued functions\\
 $F(\xi){\in}[H_{1},H_{2}],\quad |\xi|=1$ with the norm
 \be{3}{\norm{F}}_{\infty}=\mathrm{ess sup}\norm{F}<\infty,\quad |\xi|=1.\e
 We shall say that an operator $\rho{\in}[H,H]$  is \emph{strongly positive}
 if there exists  such a number $\delta>0$ that
 \be{4}\rho>{\delta}I_{H},\e
 where $I_{H}$ is the identity operator in the space $H$. The relation
 \be{5} \rho{\gg}0\e means that the operator $\rho$ is strongly positive.
 Further we use the following version of the well-known theorem (see [1] and references there).\\
 \textbf{Theorem 1.1.} \emph{Suppose given a sequence} $ \gamma_{k}{\in}[H_{1},H_{2}], \quad 1{\leq}k<{\infty}$
 \emph{and a strongly positive operator}  $\rho{\in}[H_{2},H_{2}].$  \emph{In order for there to exist an operator function} $F(\xi){\in}L_{\infty}[H_{1},H_{2}]$ \emph{such that }
 \be{6}c_{k}(F)=\frac{1}{2\pi}\int_{|\xi|=1}{\xi}^{k}F(\xi)|d{\xi}|=\gamma_{k},\quad k=1,2,...\e
 \emph{and}
 \be{7}F^{\star}(\xi)F(\xi){\leq}{\rho}^{2}\e
 \emph{it is necessary and sufficient that }
 \be{8}\Gamma^{\star}\Gamma{\leq}R^{2},\e
 \emph{where}
 \be{9}R=\mathrm{diag}\{\rho,\rho,...\}.\e
 (The integral in the right-hand side of (1.6) converges in the weak sense.)\\
 \emph{Proof.} Let us introduce the denotations
 \be{10}F_{\rho}(\xi)=F(\xi){\rho}^{-1},\quad \gamma_{k,\rho}=\gamma_{k}{\rho}^{-1}.\e
 Relations (1.6) and (1.7) take the forms
 \be{11}\frac{1}{2\pi}\int_{|\xi|=1}{\xi}^{k}F_{\rho}(\xi)|d{\xi}|=\gamma_{k,\rho},\quad k=1,2,...\e
 \emph{and}
 \be{12}F_{\rho}^{\star}(\xi)F_{\rho}(\xi){\leq}I_{H_{2}}.\e
 In case (1.11) and (1.12) the theorem is true (see [1]). Hence  in case (1.6) and (1.7) the theorem is true as well.\\
 The aim of this work is the solution of the following extremal problem.\\
 \textbf{Problem 1.2.}\emph{In the class of functions} $F(\xi){\in}[H_{1},H_{2}],\quad |\xi|=1$ \emph{satisfying condition $(1.6)$ to find the function with the least   deviation from the zero.}\\
 As a deviation measure we do not choose a number but a strictly positive matrix
 $\rho_{min}$ such that
 \be{13}F^{\star}(\xi)F(\xi){\leq}\rho_{min}^{2}.\e
 The case of the scalar matrix $\rho_{min}$ was considered in the article [1]. The transition from the scalar matrix $\rho_{min}$ to the general case considerably
 widens the class of the problems having one and only one solution. This is important both from the theoretical and the applied view points. We note that the  $\rho_{min}^{2}$ is an analogue of the eigenvalue of the operator $\Gamma^{\star}\Gamma$.
\section{Extremal problem}
In this section we consider a particular extremal problem. Namely, we try to find $\rho_{min}$ which satisfies the condition
\be{13}\Gamma^{\star}\Gamma{\leq}R_{min}^{2}, \quad R_{min}=\mathrm{diag}\{\rho_{min},\rho_{min},...\}.\e
In order to explain the notion of  $\rho_{min}$ we introduce the notations
$B_{r}=[\gamma_{2},\gamma_{3},...], \quad B_{c}=\mathrm{col}[\gamma_{2},\gamma_{3},...].$
Then the matrix $\Gamma$ has the following structure
\be{14}\Gamma=\left[\begin{array}{cc}
                      \gamma_{1} & B_{r} \\
                      B_{c} &  \Gamma_{1}
                    \end{array}\right],\e
 where
 \be{15}\Gamma_{1}=\{\gamma_{j+k}\},\quad 1{\leq}j,k<\infty,\quad \gamma_{k}{\in}[H_{1},H_{2}].\e
 It means that
  \be{16}\Gamma^{\star}\Gamma=\left[\begin{array}{cc}
                                    A_{11} &   A_{12} \\
                                      A_{12}^{\star} &  A_{22}
                                    \end{array}\right],\e
 where
 \be{17}A_{11}=\gamma_{1}^{\star}\gamma_{1}+B_{c}^{\star}B_{c}, \quad
 A_{12}=\gamma_{1}^{\star}B_{r}+B_{c}^{\star}\Gamma,\quad
 A_{22}=\Gamma_{1}^{\star}\Gamma_{1}+B_{r}^{\star}B_{r}.\e
\emph{Further we suppose that}
\be{18} R^{2}-A_{22}{\gg}0.\e
Then relation (1.8) is equivalent to the relation
\be{19}\rho^{2}{\geq}A_{11}+A_{12}(R^{2}-A_{22})^{-1}A_{12}^{\star}.\e
\textbf{Definition 2.1.} \emph{We shall call the strongly positive operator
$\rho{\in}[H_{2},H_{2}]$ a minimal solution of inequality (2.1) if the following
requirements are fulfilled:\\
1.Inequality (2.6) is valid.}\\
2.\be{20}\rho_{min}^{2}=A_{11}+A_{12}(R_{min}^{2}-A_{22})^{-1}A_{12}^{\star}.\e
It follows from (2.8) that $\rho_{min}^{2}$ coincides with the solution of the
non-linear equation
\be{21}q^{2}=A_{11}+A_{12}(Q^{2}-A_{22})^{-1}A_{12}^{\star},\e
where
$q{\in}[H_{2},H_{2}],\quad Q=\mathrm{diag}\{q,q,...\}.$
Let us note that a solution $q^{2}$ of equation (2.8) is   \emph{ an  analogue of the eigenvalue} of the operator $\Gamma^{\star}\Gamma$.\\
Now we  are giving  the method of constructing $\rho_{min}$.
 We  apply the method of successive approximation.  We put
\be{22}q_{0}^{2}=A_{11},\quad q_{n+1}^{2}=A_{11}+A_{12}(Q_{n}^{2}-A_{22})^{-1}A_{12}^{\star},\e
where
\be{23}Q_{n}=\mathrm{diag}\{q_{n},q_{n},...\},\quad n{\geq}0.\e
\emph{Further we suppose that}
\be{24}Q_{0}^{2}-A_{22}{\gg}0.\e
It follows from relations (2.10)-(2.12) that
\be{25}q_{n}^{2}{\geq}q_{0}^{2}, \quad Q_{n}^{2}{\geq}Q_{0}^{2}{\gg}0,\quad  n{\geq}0.\e
As the right-hand side of (2.10) decreases with the growth of $q_{n}^{2}$. the following assertions are true (see[7]).\\
\textbf{Lemma 2.2.}\\
1.\emph{The sequence} $q_{0}^{2}, q_{2}^{2},...$ \emph{monotonically increases  and has
the strong limit} $\underline{q}^{2}$.\\
2.\emph{The sequence} $q_{1}^{2}, q_{3}^{2},...$ \emph{monotonically decreases  and has
the strong limit} $\overline{q}^{2}$.\\
3.\emph{The inequality}
\be{26}\underline{q}^{2}{\leq}\overline{q}^{2}\e
\emph{is true.}\\
\textbf{Corollary 2.3.} \emph{If condition (2.12) is fulfilled and }
\be{27}\underline{q}^{2}=\overline{q}^{2}\e
 \emph{then}
 \be{28}\rho_{min}^{2}=\underline{q}^{2}=\overline{q}^{2}\e
 A.Ran and M.Reurings [6] investigated equation (2.10) when $A_{ij}$ are finite order matrices. Slightly changing their argumentation we shall prove that the corresponding results are true in our case as well.\\
 \textbf{Theorem 2.4.} \emph{Let} $A_{ij}$ \emph{be defined by relations (2.5) and let  condition (2.13)  be fulfilled. If the inequalities}
 \be{29}A_{11}{\geq}0,\quad A_{22}{\geq}0,\quad A_{12}A_{12}^{\star}{\gg}0 \e
 \emph{are valid, then equation (2.10) has one and only one strongly positive solution $q^{2}$ and}
 \be{30}\rho_{min}^{2}=q^{2}=\underline{q}^{2}=\overline{q}^{2}.\e
 \emph{Proof.} In view of Lemma 2.2, we have the relations
 \be{31}\underline{q}^{2}=A_{11}+A_{12}(\overline{Q}^{2}-A_{22})^{-1}A_{12}^{\star},\e
 \be{32}\overline{q}^{2}=A_{11}+A_{12}(\underline{Q}^{2}-A_{22})^{-1}A_{12}^{\star},\e
 where  $\underline{Q}=\mathrm{diag}\{\underline{q},\underline{q},...\},\quad \overline{Q}=\mathrm{diag}\{\overline{q},\overline{q},...\}.$ According to (2.14) the inequality \be{33} y=\overline{q}^{2}-\underline{q}^{2}{\geq}0\e holds.
 The direct calculation gives
 \be{34}y=B^{\star}YB,\quad Y=\mathrm{diag}\{y,y,...\},\e
 with
\be{35}B=T(I+TYT)^{-1/2}TA_{12}^{\star}\e
Here $T=(\underline{Q}_{0}^{2}-A_{22})^{-1/2}.$
Let us introduce the operator
\be{36}P=\mathrm{diag}\{p,p,...\},\quad p=\underline{q}-A_{11}.\e
From  assumption (2.12) and relation (2.17) we deduce that
\be{37}B^{\star}PB{\ll}B^{\star}(\underline{Q}^{2}-A_{22})B=p. \e
Relation (2.25) can be written in the form
\be{38}B_{1}^{\star}B_{1}{\ll}I,\quad where \quad B_{1}=P^{1/2}Bp^{-1/2}.\e
Formula (2.22) takes the form
\be{39}y_{1}=B_{1}^{\star}Y_{1}B_{1},\quad y_{1}=p^{-1/2}yp^{-1/2},\quad Y_{1}=P^{-1/2}YP^{-1/2}.\e Inequality (2.26) imply that equation (2.27) has only trivial solution
$y_{1}=0$. The theorem is proved.\\
We can omit the condition  $A_{12}A_{12}^{\star}{\gg}0,$ when
\be{40}\mathrm{dim}H_{k}=m<\infty,\quad k=1,2.\e
In this case the following assertion is true.\\
\textbf{Theorem 2.5.} \emph{Let} $A_{ij}$ \emph{be defined by relations (2.5)  and let  conditions (2.12) and (2.28) be fulfilled. If the inequalities}
 \be{41}A_{11}{\geq}0,\quad A_{22}{\geq}0, \e
 \emph{are valid, then equation (2.10) has one and only one strongly positive solution $q^{2}$ and}
 \be{42}\rho_{min}^{2}=q^{2}=\underline{q}^{2}=\overline{q}^{2}.\e
 \emph{Proof.} Let us consider the maps
 \be{43}F(q^{2})=A_{11}+A_{12}(Q^{2}-A_{22})^{-1}A_{12}^{\star},\e
 \be{44}G(q^{2})=I_{m}+U(Q^{2}-D)^{-1}U^{\star},\e
 where
 \be{45} U=q_{0}^{-1}A_{12}Q_{0}^{-1},\quad D=Q_{0}^{-1}A_{22}Q_{0}^{-1}.\e
 Fixed points $q_{F}^{2}$ and  $q_{G}^{2}$ of maps G and F respectively are connected
 by the relation
 \be{46}q_{G}=q_{0}^{-1}q_{F}q_{0}^{-1}.\e
 In view of (2.5) and (2.33) the matrix $U$ has the form $U=[u_{1},u_{2}...]$,
 where $u_{k}$ are $m{\times}m$ matrices. Let $d=\mathrm{dimker}U^{\star}.$
 The relation $x{\in}\mathrm{ker}U^{\star},\quad x{\in}{C}^{m}$ is true if and only if
 \be{47}u_{k}^{\star}x=0,\quad k=1,2,...\e
We shall use the decomposition
\be{48}((\mathrm{ker}U^{\star})^{\bot}){\bigoplus}(\mathrm{ker}U^{\star}).\e
With respect to this decomposition the matrices $u_{k}^{\star}$ and $q^{2}$ have the forms
\be{49} u_{k}^{\star}=\left[\begin{array}{cc}
                              u_{1,k}^{\star} & 0 \\
                              u_{2,k}^{\star} & 0
                            \end{array}\right],
\quad q^{2}=\left[\begin{array}{cc}
                              q_{1}^{2} & 0 \\
                              0 & I_{d}
                            \end{array}\right],\e
where  $u_{1,k}^{\star}$,  $u_{2,k}^{\star}$  and $q_{1}^{2}$  are matrices of order
$(m-d){\times}(m-d),\quad d{\times}(m-d)$ and $(m-d){\times}(m-d)$ respectively.
We note that
\be{50}q_{11}^{2}{\geq}I_{m-d}.\e
Changing the decomposition of the space ${\ell}_{2}(H_{2})$ we can represent
$U^{\star}$, $D$ and $Q^{2}$ in the forms
\be{51}U^{\star}=\left[\begin{array}{cc}
                          U_{1}^{\star}& 0 \\
                         U_{2}^{\star} & 0
                       \end{array}\right], \quad D=\left[\begin{array}{cc}
                         d_{11} & d_{12}  \\
                         d_{21} & d_{22}
                       \end{array}\right],\quad Q^{2}=\left[\begin{array}{cc}
                                                              Q_{11}^{2} & 0 \\
                                                              0 & I
                                                            \end{array}\right],\e
where  $U_{p}^{\star}=\mathrm{col}[u_{p,1}^{\star}u_{p,2}^{\star},...],\quad (p=1,2)$ and
$ Q_{11}^{2}=\mathrm{diag}\{q_{11}^{2},q_{11}^{2},...\}.$ By direct calculation we deduce that
\be{52}(Q^{2}-D)^{-1}=T\mathrm{diag}\{Q_{11}^{2}-d_{11}-d_{11}(I-d_{22})^{-1}d_{12}^{\star},I-d_{22}\}^{-1}T^{\star},
\e where
\be{53}T=\left[\begin{array}{cc}
                         I & 0  \\
                         (I-d_{22})^{-1}d_{12}^{\star} & I
                       \end{array}\right].\e
Using formulas (2.39)-(2.41)  we reduce the map $G(q^{2})$ to the form
\be{54}G_{1}(q_{11}^{2})=\hat{A}_{11}+\hat{A}_{12}
(Q_{11}^{2}-\hat{A}_{22})^{-1}\hat{A}_{12}^{\star},
\e
where
$\hat{A}_{11}=I_{m-d}+u_{2}( (I-d_{22})^{-1}u_{2}^{\star}, \quad
\hat{A}_{12}=u_{1}+u_{2}(I-d_{22})^{-1}d_{12}^{\star},$\\
$\hat{A}_{22}=d_{11}+d_{12}(I-d_{22})^{-1}d_{12}^{\star}.$
Relations (2.18), (2.33) and (2.38) imply that
\be{56}D{\ll}I,\quad Q_{11}^{2}{\geq}I.\e Hence according to (2.43) the map $G_{1}(q_{11}^{2})$ satisfies  condition (2.12). By repeating the described reduction
method we obtain the following result:
either  $\hat{A}_{12}^{\star}=0$          or $\mathrm{ker}\hat{A}_{12}^{\star}=0$
It is obvious that the theorem is true if $\hat{A}_{12}^{\star}=0$. If
 $\mathrm{ker}\hat{A}_{12}^{\star}=0$, then the $(m-d){\times}(m-d)$ matrix
 $\hat{A}_{12}\hat{A}_{12}^{\star}$ is positive, i.e. this matrix is strongly positive.
 Now the assertion of the theorem follows directly from Theorem 2.4.\\
 \textbf{Proposition 2.5.} \emph{Let conditions of either Theorem 2.3 or of Theorem 2.4 be fulfilled.Then there exists one and only one operator function
$F(\xi)$ which satisfies conditions (1.6) and (1.13)}.\\
\emph{Proof.} The formulated assertion is true when
\be{57}\rho_{min}={\alpha}I_{H_{2}},\quad \alpha=\norm{\Gamma^{\star}\Gamma}.\e
Using formulas (1.10) we reduce the general case to (2.44).The proposition is proved.\\
\textbf{Remark 2.6.}The method of constructing the corresponding operator function
is given in  paper [1] for case (2.44). Using this method we can construct the operator  function $F_{\rho_{min}}(\xi)$  and then $F(\xi)$.\\
\textbf{Remark 2.7.} Condition (2.44) is valid in a few cases. By  our approach (minimal $\rho$) we obtain the uniqueness of the solution for a broad class of problems.\\

 \end{document}